\documentclass[12pt]{article}
\usepackage{amsfonts}

\usepackage{xcolor}

\newtheorem{thm}{Theorem}[section]
\newtheorem{lemma}[thm]{Lemma}
\newtheorem{cor}[thm]{Corollary}

\newtheorem{prop}[thm]{Proposition}

\newenvironment{remark}{\par\medskip\noindent{\bf Remark.\ }}{\par\smallskip}
\newcommand{\proof
}{\par\medskip\noindent {\bf Proof.\ \ }}

\newcommand{\be}{\begin{equation}}
\newcommand{\ee}{\end{equation}}
\newcommand{\openbox}{\leavevmode
  \hbox to8pt{\hfil\vrule\vbox to6pt{\hrule width6pt\vfil\hrule}\vrule}}

\newcommand{\qed}{\hbox to5pt{ } \hfill \openbox\bigskip\medskip}

\newcommand{\rk}{\mbox{\rm rank}}
\newcommand{\ve}[1]{\mathbf{#1}}

\newcommand{\cS}{\mbox{$\cal S$}}
\newcommand{\sm}{\mbox{\rm Sm}}

\newcommand{\cB}{\mbox{$\cal B$}}

\newcommand{\cA}{\mbox{$\cal A$}}

\newcommand{\Sf}{\mathbb S}

\newcommand{\K}{\mathbb K}

\newcommand{\Z}{\mathbb Z}

\newcommand{\R}{\mathbb R}
\newcommand{\C}{\mathbb C}
\newcommand{\F}{\mathbb F}

\title{An upper bound for the size of $s$-distance sets in real algebraic sets}

\author{G\'abor Heged\"{u}s\footnote{ \'Obuda University, John von Neumann Faculty of Informatics, 
hegedus.gabor@nik.uni-obuda.hu}, 
Lajos R\'onyai\footnote{Institute of Computer Science and Control; and
Department of Algebra,
Budapest University of Technology, 
lajos@info.ilab.sztaki.hu}
}

\begin{document}
\footnotetext{Research supported in part by National Research, Development and
Innovation Office - Grant No. NKFI-115288.}

%\footnotetext{
%{\bf Keywords.} lattice polytope

%{\bf 2000 Mathematics Subject Classification.} 05A15, 13P10, 16E05 }

%\author{G\'abor Heged\"{u}s
%\\
%Lajos Rónyai 
%{\normalsize }
%}

%\begin{document}

\maketitle
\begin{abstract}

In a recent paper \cite{PP} Petrov and Pohoata developed a new algebraic method which combines 
the  Croot-Lev-Pach Lemma from additive combinatorics and Sylvester's Law of Inertia for 
real quadratic forms. As an application, they gave a simple proof  of the Bannai-Bannai-Stanton
bound on the size of $s$-distance sets (subsets $\mbox{$\cal A$}\subseteq {\mathbb R}^n$ which determine at most $s$ different distances). In this paper we extend their work and prove upper bounds for the size of 
$s$-distance sets in various real algebraic sets.  This way we obtain a novel and short proof for 
the bound of Delsarte-Goethals-Seidel on spherical $s$-distance sets and a generalization of a bound  by Bannai-Kawasaki-Nitamizu-Sato on $s$-distance sets on unions of spheres. In our arguments
we use the method of Petrov and Pohoata together with some Gr\"obner basis techniques.
\end{abstract}
\medskip



\section{Introduction}

Let  $\cA\subseteq {\R}^n$ be an arbitrary set. 
Denote by $d(\cA)$ the set of non-zero distances among the points of $\cA$:
$$
d(\cA):=\{d(\ve p_1,\ve p_2);~ \ve p_1,\ve p_2\in \cA,~ \ve p_1\ne \ve p_2\}.
$$
An {\em $s$-distance set} is a subset $\mbox{$\cal A$}\subseteq {\mathbb R}^n$ such that $|d(\mbox{$\cal A$})|\leq s$. Here we mention just two theorems from the rich area of sets with few distances, 
more information can be found for example in \cite{GY}, \cite{BB}. 
Bannai, Bannai and  Stanton proved the following upper bound for the size of an  $s$-distance set in \cite[Theorem 1]{BBS}.

\begin{thm} \label{BBSupper} 
Let $n,s\geq 1$ be integers and suppose that $\cA\subseteq {\R}^n$ is an $s$-distance set. Then
$$
|\cA|\leq {n+s\choose s}.
$$
\end{thm}

Delsarte, Goethals and  Seidel investigated $s$-distance sets on the unit sphere 
$\Sf ^{n-1}\subseteq \R^n$. These are the {\em spherical $s$-distance sets}. They  proved a general upper bound for the size of a spherical $s$-distance set in \cite{DGS}. In their proof they used Delsarte's method (see \cite[Subsection 2.2]{BB}).
\begin{thm} (Delsarte, Goethals, and  Seidel)\label{DGGupper}
Let $n,s\geq 1$ be integers and suppose that $\cA\subseteq {\Sf}^{n-1}$ is an $s$-distance set. 
Then
$$
|\cA|\leq {n+s-1\choose s}+{n+s-2\choose s-1}.
$$
\end{thm}

Before stating  our results, we introduce some notation. 
Let $\F$ be a field. In the following $S=\F[x_1, \ldots, x_n]=\F[\ve x]$ denotes  the
ring of polynomials in commuting variables $x_1, \ldots, x_n$ over $\F$.
Note that polynomials $f\in S$ can be considered as functions on $\F^n$. 
For a subset $Y$ of the polynomial ring 
$S$ and a natural number $s$ we denote by $Y_{\leq s}$ the set of polynomials from $Y$ with degree at most $s$. 
Let $I$ be an ideal of $S=\F[\ve x]$. The {\em  (affine) Hilbert
function} of the factor algebra $S/I$ is the sequence of non-negative integers $h_{S/I}(0), h_{S/I}(1),
 \ldots $, where $h_{S/I}(s)$ is the dimension over $\F$ of the factor space
 $\F[x_1,\ldots,x_n]_{\leq s}/ I_{\leq s}$ (see \cite[Section 9.3]{CLS}).
Our main technical result gives an upper bound for the size of an $s$-distance set, which is contained in a given real algebraic set.

\begin{thm}\label{maincor} Let $I\subseteq \mathbb {R}[\ve x]$ be an ideal in the polynomial ring, and 
let $\cA\subseteq \R^n $ be an $s$-distance set such that every polynomial from
$I$ vanishes on $\cA$ . Then 
$$
|\cA|\leq h_{\R[\ve x]/I}(s).
$$
\end{thm}

The proof is based on Gr\"obner basis theory and an improved version of the Croot-Pach-Lev Lemma (see \cite{CLP} Lemma 1) over the reals. Petrov and Pohoata proved this \cite[Theorem 1.2]{PP} and used it to give a new proof of Theorem \ref{BBSupper}.  We generalize their result to give a new upper bound for the size of an $s$-distance set, which is contained in a given affine algebraic set in the real affine space  $\mathbb {R}^{n}$. 

We give several corollaries, where Theorem \ref{maincor} is applied to specific 
ideals of the polynomial ring $\R[\ve x]$, the first ones being the principal ideals 
$I=(F)$, with $F\in \R[\ve x]$.

\begin{cor} \label{main3}
Let $F\in \R[\ve x]$ be a polynomial of degree $d$. Suppose that $s\geq d$. Let $\cA$ be an   $s$-distance set such that $F$ vanishes on $\cA$. Then 
$$
|\cA|\leq {n+s\choose n}-{n+s-d\choose n}.
$$
\end{cor}

For example, when $n=2$, then $F$ defines a plane curve of degree $d$. Then for $s\geq d$ we obtain  
$$
|\cA|\leq {2+s\choose 2}-{2+s-d\choose n }= 
ds-\frac{d(d-3)}{2}.$$
In particular, when $F(x,y)=y^2-f(x)$
gives a Weierstrass equation of an elliptic curve, then 
$|\cA|\leq 3s$ for 
$s\geq 3$.

\medskip

\noindent
{\bf Remark.} We can now easily derive Theorem \ref{DGGupper} for $s>1$. 
Indeed, consider the real polynomial
$$
F(x_1,\ldots ,x_n)=1-\sum_{i=1}^n x_i^2\in \R[x_1, \ldots ,x_n]
$$
of degree 2 which vanishes on ${\Sf}^{n-1}$.  Corollary \ref{main3} and the hockey-stick identity gives 
$$
|\cA|\leq {n+s \choose n}-{n+s-2 \choose n}= {n+s-1\choose s}+{n+s-2\choose s-1}.
$$

\medskip

Next, assume that  $V=\cup_{i=1}^p \cS_i$, where the  $\cS_i$ are spheres in $\R^n$.
E. Bannai, K. Kawasaki, Y. Nitamizu, and T. Sato proved the following result in \cite[Theorem 1]{BKNS} for the case when the spheres $\cS_i$ are {\em concentric}. We have a much shorter approach to the same bound, in a more general setting, without the assumption on the centers.

\begin{cor} \label{main44}
Let $\cA$ be an $s$-distance set on the union  $V$ of $p$ spheres in $\mathbb {R}^{n}$. Then
$$
|\cA|\leq \sum_{i=0}^{2p-1}  {n+s-i-1\choose s-i}. 
$$
\end{cor}

Let  $T_i\subseteq \R$ be given finite sets, where  $|T_i|=q\geq 2$ for each $i$ with  
$1\leq i\leq n$. A {\em box} is a direct product
$$
\cB:=\prod_{i=1}^n T_i\subseteq  {\R}^n.
$$
We can easily apply Theorem \ref{maincor} to obtain an upper bound for the size of $s$-distance sets in boxes.  
\begin{cor} \label{main2}
Let  $\cB\subseteq  {\R}^n$ be a box as above, and  
$\cA\subseteq \cB$  an $s$-distance set. 
Then
$$
|\cA|\leq  |\{x_1^{\alpha_1}\cdot\ldots \cdot x_n^{\alpha_n}:~ 0\leq \alpha_i \leq q-1 \mbox{ for each } i,\mbox{ and } \sum_i \alpha_i \leq s  \}|.
$$
\end{cor}

\begin{remark} In the special case $q=2$ we have  
$$
|\{x_1^{\alpha_1}\cdot\ldots \cdot x_n^{\alpha_n}:~ 0\leq \alpha_i \leq 1 \mbox{ for each } i,\mbox{ and } \sum_i \alpha_i \leq s  \}|=\sum_{j=0}^s {n\choose j},
$$
hence we obtain the upper bound 
\begin{equation} \label{binom}
|\cA|\leq \sum_{j=0}^s {n\choose j}.
\end{equation}
In the case when $T_i=T$ for $1\leq i\leq n$ and $|T|=2$, 
the Euclidean distance is essentially the same as the Hamming distance. For this case (\ref{binom}) was proved by Delsarte 
\cite{D}, see also \cite[Theorem 1]{BS}. 
\end{remark}

\begin{remark} The bound is sharp, when $q=2$, $n=2m$ and $s=m$. 
Then the 0,1 vectors of even Hamming weight give an extremal 
family $\cA\subseteq \R^n$.
\end{remark}

\begin{remark} The bound of Corollary \ref{main2}
can be nicely formulated in terms of extended binomial coefficients 
(see \cite[Example 8]{Egger} or \cite[Exercise 16]{Comtet}):
$$
|\cA|\leq \sum _{j=0}^s { n \choose  j}_q.
$$
Here ${n \choose j}_q$ is an extended binomial coefficient giving the number of restricted compositions of $j$ with $n$ terms (summands), where each term is from the set $\{0,1,\ldots ,q-1\}$.  In particular, we have 
${n \choose j}_2={n \choose j}$. 
\end{remark}

\begin{remark} 
In  \cite{H} a weaker, but similar upper bound was given for the size of $s$-distance sets in boxes:
$$
|\cA|\leq  2|\{x_1^{\alpha_1}\cdot\ldots \cdot x_n^{\alpha_n}:~ 0\leq \alpha_i \leq q-1 \mbox{ for each } i,\mbox{ and } \sum_i \alpha_i \leq s  \}|.
$$
The bound appearing in Corollary \ref{main2} presents an improvement by a factor of 2.  
\end{remark}

Let $\alpha_1,\ldots,\alpha_n$ be $n$ different elements of $\R$, and
$X_n=X_n(\alpha_1,\ldots ,\alpha_n)\subseteq \R^n$ be the set of 
permutations of $\alpha_1,\ldots ,\alpha_n$, where each permutation is considered as vector of length $n$. It was proved in \cite[Section 2]{HNR} that for $s\geq 0$
$$
h_{X_n}(s)=\sum_{i=0}^s I_n(i),
$$
where $I_n(i)$ is the number of permutations of $n$ symbols with precisely $i$ inversions. Using this, Theorem \ref{maincor} implies the following 
bound:

\begin{cor} \label{perm2}
Let $\cA\subseteq X_n$ be an  $s$-distance set. Then 
$$
|\cA|\leq \sum_{i=0}^s I_n(i).
$$
\end{cor}
\qed

In \cite[Section 5.1.1]{K} Knuth gives a generating function for $I_n(i)$ and some explicit formulae for the values $I_n(i)$, $i\leq n$.

\medskip

Let $0\leq d\leq n$ be integers and $Y_{n,d}\subseteq \R^n$
denote the set of 0,1-vectors of length $n$ which have exactly $d$ coordinate values of 1. 
The following (sharp) bound was obtained by Ray-Chaudhuri and Wilson \cite[Theorem 3]{RCW}, formulated in terms of
intersections rather than distances.

\begin{cor} \label{uniform}
Let $0\leq d\leq n$ and $s$ be integers, with  $0\leq s\leq \min(d,n-d)$. Suppose that $\cA\subseteq Y_{n,d}$ is an  $s$-distance set. Then 
$$
|\cA|\leq {n\choose s}.
$$
\end{cor}

%Let $\cP_m\subseteq \R^{m^2}$ denote the set of $m\times m$ permutation matrices.
%Corollary \ref{uniform},  applied with parameters $n=m^2$ and $d=m$ implies the 
%next bound.

%\medskip
%\begin{cor} \label{perm4}
%Suppose that $0\leq s\leq m$. Let  $\cA\subseteq \cP_m$ be an  $s$-distance set. Then 
%$$
%|\cA|\leq {m^{2}\choose  s}.
%$$
%\end{cor}
%\qed

\medskip

In some cases data about the complexification of 
a real affine algebraic set can be used to give a bound.
We give next a statement of this type. For a subset $X\subseteq \F^n$ of the affine space 
we write $I(X)$ for the ideal of all polynomials $f\in \F[\ve x]$ which vanish on $X$.

\begin{cor} \label{main_hf}
Let $V\subseteq \C^n$ be an affine variety such that the 
projective closure $\overline{V}$ of $V$ has dimension d and degree $k$.
Suppose also that the ideal $I(V)$ of $V$ is generated by polynomials over $\R$. Let
${\cal A}\subseteq V\cap \R^n $ be an $s$ distance set. Then we have
$$
|\cA|\leq \frac{k\cdot s^d}{d!}+O(s^{d-1}).
$$
\end{cor}

For instance, when in Corollary \ref{main_hf} the projective variety $\overline{V}$ is a curve of degree $k$, then the bound is 
$ks+b$ for large $s$, where $b$ is an integer. More specifically, when $\overline{V}$ is an elliptic curve such that $V\subseteq \C^2$ is the set of zeroes of 
$y^2-f(x)$, where $f(x)\in \R[x]$ is a cubic polynomial without multiple roots,
then in fact, the preceding bound becomes 
$|\cA|\leq 3s+b$ for $s$ large (see also the remark after Corollary \ref{main3}).

\medskip

The rest of the paper is organized as follows. Section 2 contains some preliminaries 
on Gr\"obner bases, Hilbert functions, and related notions. Section 3 contains the proofs of the main 
theorem and the proof of the corollaries. 

\section{Preliminaries}

%\subsection{Gr\"obner bases and standard monomials}
\label{first}

A total ordering  $\prec$ on the monomials $x_1^{i_1}x_2^{i_2}\cdots
x_n^{i_n}$  composed from
variables $x_1,x_2,\ldots, x_n$ is a {\em term order}, if 1 is the
minimal element of $\prec$, and $uw\prec vw$ holds for any monomials
$u,v,w$ with $u\prec v$. Two important  term orders are 
 the lexicographic
order $\prec_l$ and the deglex order $\prec _{dl}$. We have
$$x_1^{i_1}x_2^{i_2}\cdots x_n^{i_n}\prec_l x_1^{j_1}x_2^{j_2}\cdots
x_n^{j_n}$$
iff $i_k<j_k$ holds for the smallest index $k$ such
that $i_k\not=j_k$. As for the deglex order, we have $u\prec_{dl} v$ iff  either
$\deg u <\deg v$, or $\deg(u) =\deg(v)$, and $u\prec_lv$.

Let $\prec$ be a fixed term order. The {\em leading monomial} ${\rm lm}(f)$
of a nonzero polynomial $f$ from the ring $S=\F[\ve x]$ is the largest
(with respect to $\prec$) monomial
which occurs with nonzero coefficient in the standard form of $f$.

Let $I$ be an ideal of $S$. A finite subset $G\subseteq I$ is a {\it
Gr\"obner basis} of $I$ if for every $f\in I$ there exists a $g\in G$ such
that ${\rm lm}(g)$ divides ${\rm lm}(f)$. It can be shown that $G$ is in fact a basis of $I$.
A fundamental result is (cf. \cite[Chapter 1, Corollary
3.12]{CCS} or \cite[Corollary 1.6.5, Theorem 1.9.1]{AL}) that every
nonzero ideal $I$ of $S$ has a Gr\"obner basis with respect to $\prec$.
\medskip

A monomial $w\in S$ is a {\em standard monomial} for $I$ if
it is not a leading monomial of any $f\in I$. Let ${\rm Sm}(\prec,I,\F)$ denote
the
set of all standard monomials of $I$ with respect to the term-order 
$\prec$ over $\F$.
It is known (see \cite[Chapter 1, Section 4]{CCS}) that for a
nonzero ideal $I$  the set 
${\rm Sm}(\prec,I,\F)$ is a basis of the factor space $S/I$ over $\F$. Hence every $g\in S$ can be written uniquely as $g=h+f$ where $f\in I$ and
$h$ is a unique $\F$-linear combination of 
monomials from ${\rm Sm}(\prec,I,\F)$.

If $X\subseteq \F^n$ is a finite set, then an interpolation argument
gives that every function from $X$ to $\F$ is a polynomial function. The
latter two facts imply that

\begin{equation}\label{standard}
|\mbox{Sm}(\prec, I(X),\F)|=|X|,
\end{equation}
where $I(X)$ is the ideal of all polynomials from $S$ which vanish on $X$, and 
$\prec$ is an arbitrary term order.

The {\em initial ideal} ${ \rm in}(I)$ of $I$ is the ideal in $S$ generated
by the set of  monomials $\{ {\rm lm}(f):~f\in I\}$.

\medskip

It is easy to see \cite[Propositions 9.3.3 and 9.3.4]{CLS} that the value at $s$ of the Hilbert function $h_{S/I}$ is the 
number of standard monomials of degree at most $s$, where the ordering 
$\prec$ is deglex:
\begin{equation}
\label{hilb}
h_{S/I}(s)=|\mbox{Sm}(\prec_{dl}, I,\F)\cap \F[\ve x]_{\leq s}|.
\end{equation}
In the case when $I=I(X)$ for some  $X \subseteq \F^n$, then  
$h_{X}(s):=h_{S/I}(s)$ is the dimension of the space of functions
from
$X$ to $\F$ which are polynomials of degree at
most $s$.

\medskip

%We state here our main  combinatorial results.

\noindent

Next we recall a known fact about the Hilbert function. It concerns the change of the coefficient field. Let $\F\subset \K $ be fields and let $I\subseteq \F[\ve x]$ be an ideal, and consider the corresponding ideal $J=I\cdot \K[\ve x]$ generated by $I$ in $\K[\ve x]$. 

\begin{lemma}\label{Hf_inv} For the respective affine Hilbert functions for $s\geq 0$ we have 
$$
h_{\F[\ve x]/I}(s)=h_{\K[\ve x]/J}(s).
$$

\end{lemma}

For the convenience of the reader we outline a simple proof.

\proof
It follows from Buchberger's  criterion \cite[Theorem 2.6.6]{CLS} that a deglex Gr\"obner basis of $I$  in $\F[\ve x]$ will be a deglex Gr\"obner basis of $J$ in 
$\K[\ve x]$, implying that the initial ideals  
$\mbox{in} (I)$ and $\mbox{in}(J)$ contain exactly the same set of monomials, hence their respective factors have the same Hilbert function
$h_{\F[\ve x]/\mbox{in}(I)}(s)=
h_{\K[\ve x]/\mbox{in}(J)}(s)$, see \cite[Proposition 9.3.3]{CLS}.
Then  by \cite[Proposition 9.3.4]{CLS} we have 
$$
h_{\F[\ve x]/I}(s)=
h_{\F[\ve x]/\mbox{in}(I)}(s)=
h_{\K[\ve x]/\mbox{in}(J)}(s)=h_{\K[\ve x]/J}(s),
$$
for every integer $s\geq 0$. \qed

The projective (homogenized) version of the next statement is discussed in \cite[Example 6.10]{EH}.
\begin{prop} \label{Hyp_hilb}
Let $F\in \F[\ve x]$ be a polynomial of degree $d$. Then for  $s\geq d$ we have   
$$
h_{\F[\ve x]/(F)}(s)={n+s\choose n
}-{n+s-d\choose n}.
$$
If  $0\leq s<d$, then
$$
h_{\F[\ve x]/ (F)}(s)={n+s\choose n}.
$$
\end{prop}
\proof
By definition
$$
h_{\F[\ve x]/(F)}(s)=\dim \F[\ve x]_{\leq s}/
(F)_{\leq s}= 
$$
$$
=\dim \F[\ve x]_{\leq s}- \dim (F)_{\leq s}.
$$
Clearly
$$
\dim \F[\ve x]_{\leq s}={n+s\choose n}.
$$
Moreover
$$
(F)_{\leq s}=\{G\in \F[\ve x]_{\leq s}:~ \mbox{ there exists an  }H\in \F[\ve x] \mbox{ such that }FH=G\}.
$$
Using the fact that $\F[\ve x]$ is a domain, we see that the dimension of the latter subspace is  

$$
\dim \{H\in \R[\ve x]:~ \deg(H)\leq s-d\}=\dim \F[\ve x]_{\leq (s-d)}.
$$
The statement now follows from the fact that if $s\geq d$, then 
$$
\dim \F[\ve x]_{\leq (s-d)}={n+s-d\choose n},
$$
while for   $s<d$  we have 
$$
\dim \F[\ve x]_{\leq (s-d)}=0.
$$
\qed

\section{Proofs}

\subsection{Proof of the main result}

Petrov and Pohoata proved the following result \cite[Theorem 1.2]{PP}. They used it to give a short proof of Theorem \ref{BBSupper}. This improved version of the Croot-Lev-Pach Lemma  
has a crucial role in the proof of our results.
\begin{thm} \label{Pet} 
Let $W$ be  an $n$-dimensional vector space over a field $\F$ and let $\cA\subseteq W$ be a finite set. Let $s\geq 0$ be an integer an let $p(\ve x,\ve y)\in \F[\ve x,\ve y]$ be a $2n$-variate polynomial of degree at most $2s+1$. Consider the matrix $M(\cA,p)_{\ve a,\ve b\in A}$, where 
$$
M(\cA,p)(\ve a,\ve b)=p(\ve a,\ve b).
$$ 
This matrix corresponds to a bilinear form of ${\F}^{\cA}$ by the formula
$$
\Phi_{\cA,p}(f,g)=\sum_{\ve a,\ve b\in \cA} p(\ve a,\ve b)f(\ve a)g(\ve b), 
$$
for each $f,g:\cA\to \F$. This $\Phi_{\cA,p}$ defines a quadratic form $\Phi_{\cA,p}(f,f)$. 
In the case $\F=\R$ denote by $r_+(\cA,p)$ and  $r_-(\cA,p)$ the inertia indices of the quadratic form $\Phi_{\cA,p}(f,f)$. Then 
\begin{itemize}
\item[(i)]  $\rk(M(\cA,p))\leq 2h_{\cA}(s)$,
\item[(ii)]  if $\F=\R$, then $\max(r_+(\cA,p),r_-(\cA,p))\leq h_{\cA}(s)$.
\end{itemize}
\end{thm}

By combining Theorem  \ref{Pet} with facts about standard monomials, we have the following 
simple and elegant  upper bound for the degree of  deglex standard monomials of an  
$s$-distance set.

\begin{thm} \label{main} 
Let $\cA\subseteq \mathbb {R}^{n}$ be an $s$-distance set. Then 
$$
\mbox{Sm}(\prec_{dl}, I(\cA),\F)\subseteq 
\R[\ve x]_{\leq s}. 
$$
\end{thm}

 \noindent
 {\bf Proof.}
We follow the argument of Theorem of \cite[Theorem 1.1]{PP}.
 Let $\cA\subseteq \mathbb {R}^{n}$ denote an $s$-distance set. Recall that $d(\cA)$ denotes the set of (non-zero) distances among points of $\cA$. Define the 
$2n$--variate polynomial by:
$$
p(\ve x,\ve y)=\prod_{t\in d(\cA)} \Big(t^2-\Vert \ve x-\ve y \Vert^2 \Big)\in \R[\ve x,\ve y].
$$
Then
we can apply Theorem  \ref{Pet} for $p(\ve x,\ve y)$ whose degree is $2s$. The matrix $M(\cA,p)$ is a positive diagonal matrix, giving that 
$$
r_+(\cA,p)=|\cA|.
$$
It follows from Theorem  \ref{Pet} (ii) that 
$$
|\cA|=r_+(\cA,p)\leq h_{\cA}(s).
$$
But  equations (\ref{hilb}),
(\ref{standard}) and the finiteness of $\cA$ imply that 
$$
|\cA|\leq h_{\cA}(s)=|\sm(\prec_{dl},I(\cA),\R)\cap
\R[\ve x]_{\leq s}|\leq |\sm(\prec_{dl},I(\cA),\R)|=|\cA|.
$$
We infer that  
$$
 |\sm(\prec_{dl},I(\cA),\R)\cap\R[\ve x]_{\leq s}|= |\sm(\prec_{dl},I(\cA),\R)|,
$$
and hence 

$$
\mbox{Sm}(\prec_{dl}, I(\cA),\R)\subseteq 
\R[\ve x ] _{\leq s}. 
$$
\qed

\noindent
{\bf Proof of  Theorem  \ref{maincor}.}
Theorem \ref{main} gives that
$$
\mbox{Sm}(\prec_{dl}, I(\cA),\R)\subseteq 
\R[\ve x]_{\leq s}. 
$$
Since $I$ vanishes on $\cA$, we have $I\subseteq I(\cA)$, hence 
$$
\mbox{Sm}(\prec_{dl}, I(\cA),\R)\subseteq \mbox{Sm}(\prec_{dl}, I,\R).
$$
The preceding two equations imply that 
$$
\mbox{Sm}(\prec_{dl}, I(\cA),\R)\subseteq \mbox{Sm}(\prec_{dl}, I,\R)\cap \R[\ve x]_{\leq s}. 
$$
Now it follows from (\ref{hilb}) and 
(\ref{standard}) that
$$
|\cA|=|\mbox{Sm}(\prec_{dl}, I(\cA),\R)|\leq |\mbox{Sm}(\prec_{dl}, I,\R)\cap \R[\ve x]_{\leq s}|= h_{\R[\ve x]/I}(s).
$$
\qed

\subsection{Proofs for the Corollaries}

\noindent 
{\bf Proof of Corollary \ref{main3}}.
From Theorem \ref{maincor} we obtain the bound 
$ |\cA|\leq h_{\R[\ve x]/(F)}(s)$, therefore 
for  $s\geq d$ we have 
$$
|\cA|\leq h_{\R[\ve x]/(F) }(s)={n+s\choose n}-{n+s-d\choose n},
$$
by Proposition \ref{Hyp_hilb}. \qed

\noindent
{\bf Proof of Corollary \ref{main44}}.
It is easy to verify that 
$$
 \sum_{i=0}^{2p-1}  {n+s-i-1\choose s-i}={n+s\choose s}-{n+s-2p\choose n}.
$$

Let $V=\cup_{i=1}^p \cS_i$, and
assume, that the center of the sphere $\cS_i$ is the point 
$(a_{1,i},\ldots ,a_{n,i})\in \mathbb {R}^{n}$ and the radius of $\cS_i$ is $r_i\in \R$ 
for  $i=1,\ldots , p$.  
Next consider the polynomials 
$$
F_i(x_1, \ldots, x_n)=(\sum_{m=1}^n (x_m-a_{m,i})^2 )-r_i^2\in \R[x_1, \ldots ,x_n]
$$
for each $i$ and put $F:=\prod_i F_i$. Then $\deg(F)=2p$ and  $F$ vanishes on $V$. We  may apply Corollary \ref{main3} for the polynomial $F$. Then for 
$s\geq 2p$ we obtain the desired bound
$$
|\cA|\leq {n+s\choose n}-{n+s-2p\choose n}.
$$
When $s<2p$, the bound follows from the Bannai-Bannai-Stanton theorem. 
\qed

\noindent 
{\bf Proof of  Corollary \ref{main2}:}
It is  well-known and easily proved that the following set of polynomials is a (reduced) Gr\"obner basis of the ideal $I(\cB)$  (with respect to any term order):
$$
\{\prod_{t\in T_i} (x_i-t):~ 1\leq i\leq n\}.
$$
This readily gives the (deglex) standard monomials for $I(\cB)$: 
$$
\mbox{Sm}(\prec_{dl}, I(\cB),\R)=|\{x_1^{\alpha_1}\cdot\ldots \cdot x_n^{\alpha_n}:~ 0\leq \alpha_i \leq q-1 \mbox{ for each } i \}|.
$$
It follows from  Theorem  \ref{maincor}  and equation (\ref{hilb}) that
$$
|\cA|\leq h_{\cB}(s)=|\mbox{Sm}(\prec_{dl}, I(\cB),\R)\cap \R[\ve x]_{\leq s}|=
$$
$$
=|\{x_1^{\alpha_1}\cdot\ldots \cdot x_n^{\alpha_n}:~ 0\leq \alpha_i \leq q-1 \mbox{ for each } i,\ \mbox{ and } \sum_i \alpha_i \leq s  \}|.
$$
\qed

\noindent
{\bf Proof of Corollary \ref{uniform}.}
The statement follows at once from the result 
\begin{equation}
h_{Y_{n,d}}(s)={n \choose s}.
\end{equation}
proved by Wilson in \cite{W} (formulated there in the language of inclusion matrices, see also \cite[Corollary 3.1]{HR}),
and Theorem \ref{maincor}. 
\qed
\medskip

\noindent 
{\bf Proof of Corollary \ref{main_hf}}.
Write $I=I(V)\cap \R[\ve x]$ and $J =I(V)\subseteq \C[\ve x]$.
It follows from Theorem \ref{maincor} and Proposition \ref{Hf_inv}  that 
$$
|\cA|\leq h_{\R[\ve x]/I}(s)= h_{\C[\ve x]/J}(s).
$$

From Theorem 9.3.12 of \cite{CLS} we obtain that 
the affine Hilbert function $h_{\C[\ve x]/J}(s) $ is the same as the projective Hilbert function 
$h_{\overline{V}}(s)$ of the projective variety 
$\overline{V}$. 
Now \cite[Proposition 13.2]{H2} and the subsequent remark imply that for $s$ large the Hilbert function will be the same as the Hilbert polynomial:  $h_{\overline{V}}(s) =p_{\overline{V}}(s)$, moreover 
$$
 p_{\overline{V}}(s)=\frac{k}{d!}\cdot s^d+\mbox{ terms  of degree at most }  d-1 \mbox{ in } s.
$$
This proves the statement.
\qed

%\subsection{Gr\"obner  theory}

% {\bf Acknowledgements.} 

\end{document}